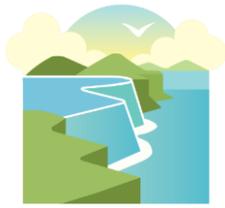
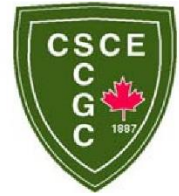

# MODELLING OF UNSATURATED FLOW THROUGH POROUS MEDIA USING MESHLESS METHODS

Boujoudar, M.[1,4], Beljadid, A.[1,2,5], and Taik, A.[3,6]
[1] Mohammed VI Polytechnic University, Green City, Morocco
[2] University of Ottawa, Ottawa, Canada
[3] FST-Mohammedia, University Hassan II, Casablanca, Morocco
[4] mohamed.boujoudar@um6p.ma
[5] abeljadi@uttawa.ca
[6] ahmed.taik@fstm.ac.ma

**Abstract:** In this study, we focus on the modelling of infiltration process in porous media. We use the meshless techniques for efficiently solving the Richards equation which describes unsaturated water flow through soils. The design of approximate numerical methods for the Richards equation remains computationally challenging and requires the development of efficient numerical techniques. This difficulty is mainly due to the nonlinearity of the unsaturated hydraulic conductivity and the capillary pressure function. In this study, we develop a new method based on the localized radial basis function (RBF) and the Kirchhoff transformation technique in order to solve Richards equation in one and two-dimensional homogeneous medium. Our approach using the multiquadric radial basis function allows us to reduce the computational time and provide accurate numerical solutions. The proposed method does not require mesh generation. Picard's iterations are used to linearize the resulting nonlinear problem obtained using the Kirchhoff transformation technique. The numerical simulations show the capability of the proposed numerical techniques in predicting the dynamics of water through unsaturated soils.

## 1   INTRODUCTION

The process of infiltration through porous media is an important part of hydrological cycle. The modelling of this process has important practical applications in engineering such as water resources management and agriculture. Most numerical models that describe the unsaturated flow in soils use the Richards model (Richards 1931) which is a highly nonlinear equation. This equation is obtained from Darcy's law and the conservation of mass (Bear 1972). The strong non-linearity of the unsaturated conductivity and the capillary pressure as functions of saturation and the presence of both advection and diffusion terms make the Richards equation more challenging in terms of numerical approximations and require the development of efficient numerical techniques. The unsaturated conductivity and capillary pressure are correlated using empirical models and experiment data such as the van Genuchten (Van Genuchten 1980), Brooks-Corey (Brooks and Corey 1964) and Gardner (Gardner 1958) models. Analytical solutions of the Richards equation can only obtained for some cases with special initial and boundary conditions (Sander, et al. 1988, Srivastava and Yeh 1991, Tracy 2011, Huang and Wu 2012). Therefore, different numerical techniques are developed to efficiently solve the Richards equation such as, finite-difference, finite-element, and finite-volume methods. For instance, (Celia, Bouloutas and Zarba 1990) used the mixed form of the Richards equation and proposed a general mass-conservative numerical scheme, and (Bause and Knabner 2004) developed an adaptive mixed hybrid finite element discretization for the Richards equation. (Manzini and Ferraris 1990) developed a mass conservative finite volume method using two-dimensional unstructured grids. Although many numerical techniques have been developed to numerically solve the Richards equation, there is still a strong need for more robust numerical techniques for modelling flows in unsaturated soils.

The aim of this work is to develop a new technique based on the localized radial basis function method and the Kirchhoff transformation in order to solve Richards equation in one and two-dimensional homogeneous medium. The proposed technique allows us to avoid mesh generation, which makes the numerical method less expensive in terms of computational cost. The use of localized meshless method has the advantage of flexibility in dealing with complex geometries (Boujoudar, Beljadid and Taik 2021). The proposed method performs well in terms of accuracy and efficiency for modelling unsaturated flow through soils.

To handle the nonlinearity of the Richards equation, we use the Kirchhoff transformation which allows us to reduce the nonlinearity of the studied problem. We used Picard iterations to solve the problem with the Kirchhoff variable where we used the backward Euler method for temporal discretization. Other numerical techniques using the Kirchhoff transformation to solve the Richards equation can be found in (Ji, et al. 2008). The performance of the proposed numerical method is assessed using different test cases.

The outline of the paper is as follows. In section 2, we introduce the governing equation and the proposed system using the Kirchhoff transformation. In section 3, we present the proposed meshless method. Numerical simulations are performed in section 4 for modelling water flow through one and two-dimensional unsaturated porous media. Finally, we provide some conclusions in section 5.

## 2   GOVERNING EQUATION

### 2.1   The mathematical model

Infiltration of water in unsaturated soils is described by the Richards equation (Richards 1931) which can be derived from Darcy's law and the conservation of mass. This equation is given by:

$$[1] \quad \frac{\partial \theta}{\partial t} + \nabla \cdot (K \nabla h) + \frac{\partial K}{\partial z} = s(\boldsymbol{x}, t), \quad \boldsymbol{x} \in \Omega, 0 \leq t \leq T,$$

where $\theta$ $[L^3/L^3]$ is the moisture content, $h$ [L] is the pressure head, $K$ $[L/T]$ is the unsaturated hydraulic conductivity, $\boldsymbol{x} = (x, y, z)^T$ is the coordinate vector, $x$ [L] and $y$ [L] denote the horizontal dimensions and $z$ [L] denotes the vertical dimension positive down (coordinate in the direction of gravity) and $s(\boldsymbol{x}, t)$ is a source or sink term which can depend on space and time and can include evaporation, plant root extraction, etc. In this study, we assume that $s(\boldsymbol{x}, t) = 0$, $\Omega$ is an open set of $\mathbb{R}^d (d = 1,2,3)$, and $T$ is the final simulation time.

We note that the Richards equation can be expressed using the water saturation $S = \left(\frac{\theta - \theta_r}{\theta_s - \theta_r}\right)$ and the parameter $\phi = \theta_s - \theta_r$ where $\theta_s$ and $\theta_r$ are respectively the saturated and residual moisture contents. The unsaturated hydraulic conductivity is given by:

$$[2] \quad K = K_s k_r,$$

where $k_r$ is the water relative permeability, which accounts for the effect of partial saturation and the saturated hydraulic conductivity is as follows:

$$[3] \quad K_s = \frac{\rho g k}{\mu},$$

where $\rho$ is the water density, $g$ is the gravitational acceleration, $k$ is the intrinsic permeability of the medium, and $\mu$ is the fluid dynamic viscosity. The Richards equation can be rewritten in the following form:

$$[4]\ \phi \frac{\partial S}{\partial t} + \nabla \cdot (K_s k_r \nabla h) + \frac{\partial (K_s k_r)}{\partial z} = 0, \quad x \in \Omega, 0 \leq t \leq T,$$

Eq. 4 is highly non-linear due to the nonlinearity of the hydraulic conductivity and the capillary pressure function. Constitutive relationships are available for the functions $S\ [L^3/_{L^3}]$ and $K\ [L/_T]$ based on experiment. In our study, the numerical techniques will be developed based on Eq. 4 where we will introduce the Kirchhoff transformation in order to reduce the nonlinearity of the equation.

## 2.2 Capillary pressure

The pressure head can be expressed as a function of saturation in the following form:

$$[5]\ h(S) = h_{cap} J(S),$$

where $J(S)\ [-]$ is a dimensionless capillary pressure function and $h_{cap}\ [L]$ is the capillary rise which is given by the classical Leverett scaling (Leverett 1941):

$$[6]\ h_{cap} \sim \frac{\gamma \cos \boldsymbol{\theta}}{pg \sqrt{\frac{k}{\phi}}},$$

$\gamma$ is the surface tension between the fluids, and $\boldsymbol{\theta}$ is the contact angle.

## 2.3 Kirchhoff transformation

The Kirchhoff integral transformation is defined as:

$$[7]\ \varphi(h) = \int_{+\infty}^{h} k_r(s)\, ds,$$

and by applying this transformation, we can rewrite the Richards equation using the variable $\varphi$, as explained below:

$$[8]\ \nabla \cdot (K_s k_r \nabla h) = K_s \nabla^2 \varphi,$$

$$[9]\ \frac{\partial \varphi}{\partial t} = k_r \frac{\partial h}{\partial t}, \quad \frac{\partial \varphi}{\partial z} = k_r \frac{\partial h}{\partial z}.$$

By transforming the derivative terms $\frac{\partial S}{\partial t}$ and $\frac{\partial}{\partial z}(K_s k_r)$ using the variable $\varphi$, we obtain:

$$[10]\ \frac{\partial S}{\partial t} = \frac{\partial S}{\partial h}\frac{\partial h}{\partial t} = \left(k_r^{-1} \frac{\partial S}{\partial h}\right) \frac{\partial \varphi}{\partial t},$$

$$[11]\ \frac{\partial}{\partial z}(K_s k_r) = K_s \frac{\partial k_r}{\partial h}\frac{\partial h}{\partial z} = \left(K_s k_r^{-1} \frac{\partial K_r}{\partial h}\right) \frac{\partial \varphi}{\partial z}.$$

We consider the variables:

$$[12]\ \begin{cases} A = -\dfrac{\phi}{K_s}\left(k_r^{-1} \dfrac{\partial S}{\partial h}\right), \\ B = \left(K_r^{-1} \dfrac{\partial K_r}{\partial h}\right), \end{cases}$$

this leads to the following equation:

[13] $A \frac{\partial \varphi}{\partial t} + \nabla^2 \varphi + B \frac{\partial \varphi}{\partial z} = 0.$

Finally, by applying the Kirchhoff transformation, we reduced the nonlinearity of Eq. 4 and obtain Eq. 13 which has many benefits in terms of convergence of the proposed numerical method.

### 2.4 Initial and boundary conditions

For the initial condition, we assume that the pressure head is $h(x, 0) = h_0$ for each point $x$ on the computational domain $\Omega$, which can be expressed using the Kirchhoff variable as $\varphi(x, 0) = \varphi(h_0)$.

We transform the boundary conditions using the Kirchhoff variable in a similar way:

Dirichlet: $h(x, t) = h_d$ for each $x \in \partial \Omega$ leads to $\varphi(x, t) = \varphi(h_d)$.

Neumann: $n_i \frac{\partial h}{\partial x_i} = h_N$ implies $n_i \frac{\partial \varphi}{\partial x_i} = K_r h_N$, with $h_N$ is a given function and $n_i$ is the unit normal.

## 3 THE MATERIALS AND PROPOSED TECHNIQUES

In this section, we propose an efficient computational technique based on radial basis function collocation method (Kansa 1990, Kansa1 1990). This method has recently become very popular due to its advantages in terms of approximation properties of solutions and its less computational cost since it does not require mesh generation.

Eq. 13 is solved using the localized RBF collocation method and the Picard iteration technique. The temporal discretization of Eq. 13 using the backward Euler method is given by:

[14] $A^{n+1} \frac{\varphi^{n+1} - \varphi^n}{\Delta t} + \nabla^2 \varphi^{n+1} + B^{n+1} \frac{\partial \varphi^{n+1}}{\partial z} = 0,$

where $\varphi^{n+1}$, $A^{n+1}$ and $B^{n+1}$ are the approximate values of $\varphi$, $A$ and $B$ at $t = t^{n+1}$, respectively. $\Delta t = t^{n+1} - t^n$ is the time setup and the solution is assumed to be known at $t^n$ and unknown at $t^{n+1}$.

Eq. 14 is linearized using the Picard iteration method which involves sequential estimation of the unknown $\varphi^{n+1}$ using the latest estimates of $A^{n+1}$ and $B^{n+1}$. If $m$ identifies iteration levels, then the Picard iteration steps can be written as:

[15] $A^{m,n+1} \frac{\varphi^{m+1,n+1} - \varphi^n}{\Delta t} + \nabla^2 \varphi^{m+1,n+1} + B^{m,n+1} \frac{\partial \varphi^{m+1,n+1}}{\partial z} = 0.$

For the sake of simplicity, we consider the following notations:

[16] $\mathcal{L}^m = \left( \frac{A^{m,n+1}}{\Delta t} \cdot + \nabla^2 \cdot + B^{m,n+1} \frac{\partial \cdot}{\partial z} \right),$

[17] $f^{m,n+1} = A^{m,n+1} \frac{\varphi^n}{\Delta t},$

$\mathcal{L}^m$ is a linear operator for each Picard iteration $m$. Subject to boundary and initial conditions, Eq. 15 can be rewritten in the following form:

[18] $\begin{cases} \mathcal{L}^m \varphi^{m+1,n+1}(x) = f^{m,n+1}(x), & x \in \Omega, \\ \mathcal{B}\varphi^{m+1,n+1}(x) = q(x), & x \in \partial \Omega, \\ \varphi^{m+1,0}(x) = \varphi_0^{m+1}(x), & x \in \Omega, \end{cases}$

$\mathcal{B}$ is a border operator, $q$ is the given function associated with the boundary conditions. For each iteration $n$, Eq. 18 is solved using localized RBF meshless method at each Picard iteration $m$ until the stop condition is verified which is given by:

$$[19] \quad \delta^m = |\varphi^{m+1,n+1} - \varphi^{m,n+1}| \leq Tol,$$

with $Tol$ is the error tolerance.

### 3.1 Localized RBF meshless method

In this section, we present the local multiquadric (LMQ) method (Lee, Liu and Fan 2003). This approach is different from the traditional global multiquadric approximation since only local configuration of nodes are used. To recall the localized RBF techniques, let $\{x_j\}_{j=1}^{n_i}$ and $\{x_j\}_{j=n_i+1}^{N}$ be the collocation points in $\Omega$ and $\partial\Omega$, respectively. $n_i$ is the number of interior points and $N$ the total number of collocation points distributed over the computational domain. For each $x_s \in \Omega$, we create a localized domain $\Omega^{[s]}$ that contains $n_s$ nearest neighbors interpolation points $\{x_k^{[s]}\}_{k=1}^{n_s}$ to $x_s$.

In each localized domain $\Omega^{[s]}$, the approximate solution can be written as a linear combination of $n_s$ multiquadric functions in the following form:

$$[20] \quad \varphi_{[s]}^{m+1,n+1}(x_s) = \sum_{k=1}^{n_s} \alpha_k^{m+1,n+1} \Phi_k\left(\left\|x_s - x_k^{[s]}\right\|\right),$$

where $\{\alpha_k^{m+1,n+1}\}_{k=1}^{n_s}$ unknown coefficients to be determined, $\|.\|$ is the Euclidian norm and $\Phi_k$ are the multiquadric radial basis functions defined as:

$$[21] \quad \Phi_k(x) = \Phi(r_k) = \sqrt{1 + (\varepsilon r_k)^2},$$

where $\varepsilon > 0$ is the shape parameter and $r_k = \|x - x_k\|$. Eq. 21 can be presented in the matrix form:

$$[22] \quad \varphi_{[s]}^{m+1,n+1} = \Phi^{[s]} \alpha_{[s]}^{m+1,n+1},$$

where $\varphi_{[s]}^{m+1,n+1} = \left[\varphi_{[s]}^{m+1,n+1}(x_1^{[s]}), \varphi_{[s]}^{m+1,n+1}(x_2^{[s]}), \ldots, \varphi_{[s]}^{m+1,n+1}(x_{n_s}^{[s]})\right]^T$, $\alpha_{[s]}^{m+1,n+1} = \left[\alpha_{[s]}^{m+1,n+1}(x_1^{[s]}), \alpha_{[s]}^{m+1,n+1}(x_2^{[s]}), \ldots, \alpha_{[s]}^{m+1,n+1}(x_{n_s}^{[s]})\right]^T$ and $\Phi^{[s]}$ is an $n_s \times n_s$ real symmetric coefficient matrix defined as $\Phi^{[s]} = \left[\Phi(\|x_i^{[s]} - x_j^{[s]}\|)\right]_{1 \leq i,j \leq n_s}$. The vector $\alpha_{[s]}^{m+1,n+1}$ can be obtained as the following equation:

$$[23] \quad \alpha_{[s]}^{m+1,n+1} = \left(\Phi^{[s]}\right)^{-1} \varphi_{[s]}^{m+1,n+1},$$

for $x_s \in \Omega$, we apply the differential operator $\mathcal{L}^m$ to Eq. 20 to obtain the following equation:

$$[24] \quad \mathcal{L}^m \varphi_{[s]}^{m+1,n+1}(x_s) = \sum_{k=1}^{n_s} \alpha_k^{m+1,n+1} \mathcal{L}^m \Phi_k\left(\left\|x_s - x_k^{[s]}\right\|\right) = \sum_{k=1}^{n_s} \alpha_k^{m+1,n+1} \Psi^m\left(\left\|x_s - x_k^{[s]}\right\|\right),$$
$$= \Theta_{[s]}^m \alpha_{[s]}^{m+1,n+1} = \Theta_{[s]}^m \left(\Phi^{[s]}\right)^{-1} \varphi_{[s]}^{m+1,n+1} = \Lambda_{[s]}^m \varphi_{[s]}^{m+1,n+1} = \Lambda^m \varphi^{m+1,n+1},$$

where $\varphi^{m+1,n+1} = [\varphi^{m+1,n+1}(x_1), \varphi^{m+1,n+1}(x_2), \ldots, \varphi^{m+1,n+1}(x_N)]^T$, $\Theta_{[s]}^m = \left[\Psi^m\left(\left\|x_s - x_1^{[s]}\right\|\right), \Psi^m\left(\left\|x_s - x_2^{[s]}\right\|\right), \ldots, \Psi^m\left(\left\|x_s - x_{n_s}^{[s]}\right\|\right)\right]^T$, $\Psi^m = \mathcal{L}^m \Phi_k$ and $\Lambda_{[s]}^m = \Theta_{[s]}^m \left(\Phi^{[s]}\right)^{-1}$.

In order to extend Eq. 24 to be able to use $\varphi^{m+1,n+1}$ instead of $\varphi_{[s]}^{m+1,n+1}$, we consider $\Lambda^m$ as the expansion of $\Lambda_{[s]}^m$ which can be obtained by padding the local vector with zeros.

Similarly, for $x_s \in \partial\Omega$, we create an influence domain $\Omega^{[s]}$ containing $x_s$. Then we have:

$$[25] \mathcal{B}\varphi^{m+1,n+1}(x_s) = \sum_{k=1}^{n_s} \alpha_k^{m+1,n+1} \mathcal{B}\Phi_k\left(\left\|x_s - x_k^{[s]}\right\|\right) = \left(\mathcal{B}\Phi^{[s]}\right)\alpha_{[s]}^{m+1,n+1},$$
$$= \left(\mathcal{B}\Phi^{[s]}\right)\left(\Phi^{[s]}\right)^{-1} \varphi_{[s]}^{m+1,n+1} = \sigma^{[s]} \varphi_{[s]}^{m+1,n+1} = \sigma\varphi^{m+1,n+1},$$

where $\sigma^{[s]} = \left(\mathcal{B}\Phi^{[s]}\right)\left(\Phi^{[s]}\right)^{-1}$, and $\sigma$ is the expansion of $\sigma^{[s]}$ obtained by completing the local vector with zeros.

We substitute the equations 24 and 25 into Eq. 18 to obtain the following system:

$$[26] \begin{cases} \mathcal{L}^m \varphi^{m+1,n+1}(x_s) = \Lambda^m(x_s)\varphi^{m+1,n+1} = f^{m+1,n+1}(x_s), \\ \mathcal{B}\varphi^{m+1,n+1}(x_s) = \sigma(x_s)\varphi^{m+1,n+1} = q(x_s), \end{cases}$$

which leads to the following sparse system:

$$[27] \begin{pmatrix} \Lambda^m(x_1) \\ \Lambda^m(x_2) \\ \vdots \\ \Lambda^m(x_{n_i}) \\ \sigma(x_{n_i+1}) \\ \vdots \\ \sigma(x_N) \end{pmatrix} \begin{pmatrix} \varphi^{m+1,n+1}(x_1) \\ \varphi^{m+1,n+1}(x_2) \\ \vdots \\ \varphi^{m+1,n+1}(x_{n_i}) \\ \varphi^{m+1,n+1}(x_{n_i+1}) \\ \vdots \\ \varphi^{m+1,n+1}(x_N) \end{pmatrix} \begin{pmatrix} f^{m+1,n+1}(x_1) \\ f^{m+1,n+1}(x_2) \\ \vdots \\ f^{m+1,n+1}(x_{n_i}) \\ q(x_{n_i+1}) \\ \vdots \\ q(x_N) \end{pmatrix}.$$

The matrix generated by the localized RBF is sparse due to the presence of the local configuration in the solution approximation. This allows us to avoid ill-conditioned issues that arise in dense systems of equations generated by the global approach. By solving Eq. 27, we obtain the approximate values of $\varphi^{m+1,n+1}$ at all nodes $\varphi^{m+1,n+1}(x_s)$, $s = 1,2,\ldots,N$ of the computational domain.

## 4 NUMERICAL TESTS

In this section, we perform numerical experiments for solving the Richards equation by using the obtained Eq. 27 in one and two-dimensional systems. We used the localized RBF method based on the multiquadric radial basis functions. For the temporal discretization, we used the backward Euler method.

### 4.1 One dimensional infiltration problem

In this numerical test, we used the Brooks Corey model (Brooks and Corey 1964) which describes the pressure head and the power law for the relative permeability.

$$[28] \quad J(S) = S^{-1/\lambda}, \qquad k_r = S^\beta = \begin{cases} \left(\dfrac{h}{h_{cap}}\right)^{-\lambda\beta}, & \text{if } h \geq h_{cap}, \\ 1, & \text{if } h < h_{cap}. \end{cases}$$

Where $\lambda$ and $\beta$ are respectively the parameters related to the Brook Corey model and the power law for the relative permeability. The second (inequality) condition for the capillary pressure in Eq. 28 is introduced to avoid numerical issues (Ma, Hook et Ahuja 1999, Lenhard et Parker 1989).

We consider different types of soils with a depth $L$ and their hydraulic parameters are shown in Table 1. We simulate a one-dimensional infiltration problem using the proposed method. In order to verify the effectiveness of the developed numerical model, we compare our numerical results with the numerical solutions of 1D-Hydrus where we consider the following initially and boundary conditions:

$$[29] \begin{cases} \theta(z,0) = \theta_0, \\ \theta(0,t) = \theta_s, \\ \theta(L,t) = \theta_0, \end{cases}$$

Table 1: Hydraulic property parameters of 2 types of soil.

| Texture | $\theta_r$ | $\theta_s$ | $\theta_0$ | $K_s$ | $h_{cap}$ | $\lambda$ | $\lambda\beta$ |
|---|---|---|---|---|---|---|---|
| Sandy clay | 0.109 | 0.321 | 0.121 | 0.002 | 29.15 | 0.168 | 2.504 |
| Loam | 0.027 | 0.463 | 0.040 | 0.022 | 11.15 | 0.220 | 2.660 |

Figure 1 shows the numerical solutions obtained using the proposed method and 1D-Hydrus solutions where we observe good agreement between the solutions.

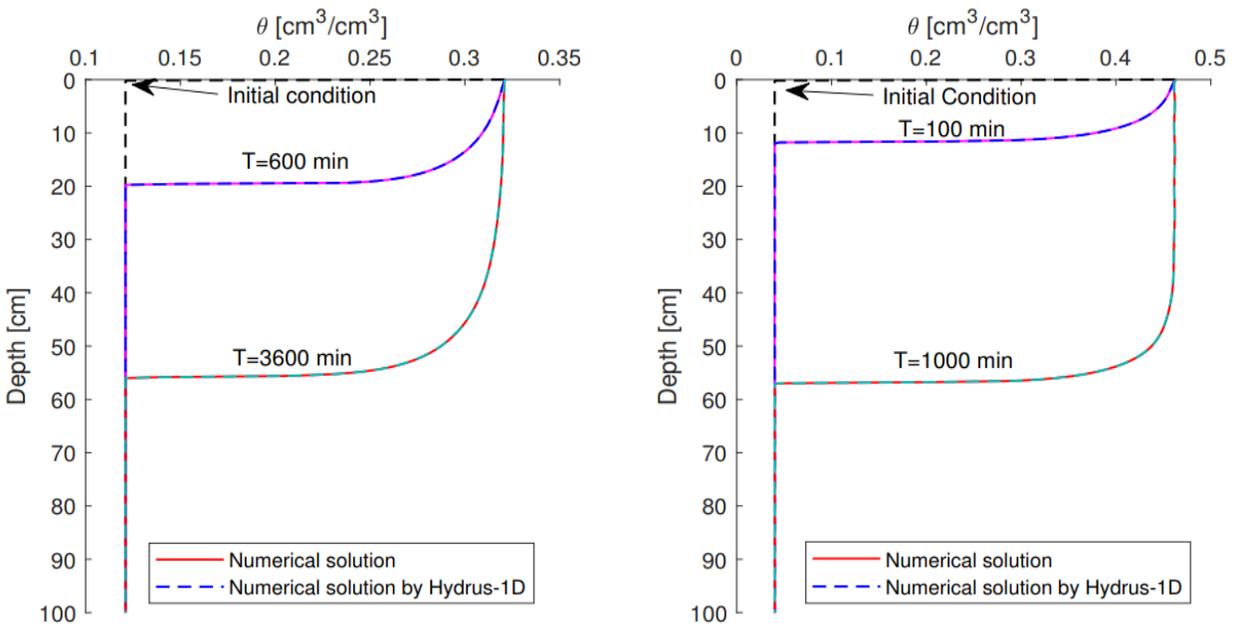

Figure 1: Time evolution of moisture content for the soils given in Table 1. Left (Sandy clay), right (Loam).

Table 2 presents the root mean squared error (RMSE), the relative error $(L_{er}^1)$ between the numerical solutions and the solutions simulated by 1D-Hydrus. The results confirm the accuracy of the proposed method in modelling unsaturated flow in soils.

Table 2: The RMSE and the $L^1_{er}$ between the numerical solutions and 1D-Hydrus solutions.

| Soils | $T$ (min) | RMSE | $L^1_{er}$ |
|---|---|---|---|
| Sandy Clay | 600 | $5 \times 10^{-3}$ | $3.5 \times 10^{-3}$ |
| | 3600 | $5.8 \times 10^{-3}$ | $4.3 \times 10^{-3}$ |
| Loam | 100 | $4.8 \times 10^{-3}$ | $1.08 \times 10^{-3}$ |
| | 1000 | $6 \times 10^{-3}$ | $7.2 \times 10^{-3}$ |

## 4.2 Two-dimensional infiltration problem

In this example, we perform numerical simulations using the proposed method for a two-dimensional infiltration problem where we consider a rectangular domain $[0, l] \times [-L, 0]$. We used the same hydraulic parameters of test 1 (Table 1) and $l = L = 100$ cm. We consider the following initial and boundary conditions:

[30] $\begin{cases} \theta(x, z, 0) = \theta_0, \\ \theta(x, 0, t) = \theta_s, \\ \theta(x, L, 0) = \theta_0, \end{cases}$

and no-flux boundary conditions are imposed on the sides $x = 0$ and $x = l$ of the domain.

The numerical simulations are performed using $N_x = 200$, $N_z = 200$, $\Delta t = 0.05$, and the localized RBF parameters $\varepsilon = 0.6$ and $n_s = 5$. The sandy clay and loam soils are selected in this numerical test to simulate unsaturated flow through a two-dimensional homogeneous medium. The time evolution of the total mass per unit of length of the 2D numerical solutions and the solutions simulated by 1D-Hydrus for a computational domain of unit length (1D problem) are displayed in Figure 2. We observe good agreement between the solutions, which demonstrates the accuracy of the proposed method in modelling 2D unsaturated flow in soils. Figure 3 shows the time evolution of saturation for the sandy clay and loam considered soils.

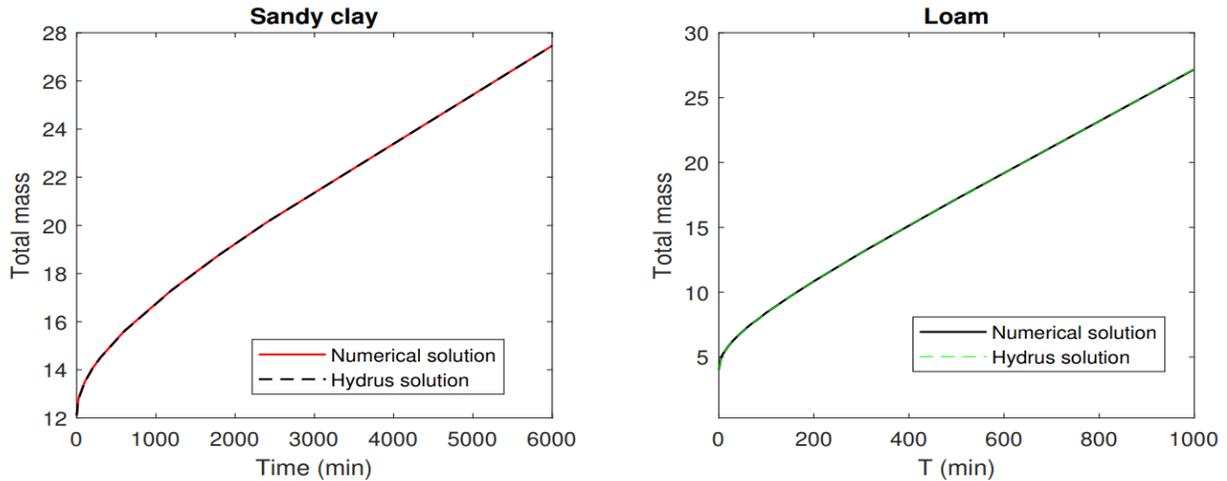

Figure 2: Time-evolution of the total mass per unit of length of the sandy clay and loam soils.

The proposed method is efficient and accurate for solving the Richards equation. The method can be used for modelling unsaturated flow through homogeneous soils.

## 5 Conclusion

This paper focused on the infiltration process in porous media and introduced computational techniques for efficiently solving the Richards equation in one- and two-dimensional homogeneous medium. The proposed techniques using the Kirchhoff transformation allow us to reduce the nonlinearity of the obtained system from the Richards equation. Our approach using a localized radial basis function method avoiding mesh generation allows us to reduce the computational cost. The accuracy of the proposed method was validated using comparison between the numerical solutions and the results of 1D-Hydrus. Our results confirm the accuracy of the proposed techniques and their efficiency in terms of computational cost for solving the Richards equation. The numerical techniques proposed in this study for modelling unsaturated flow through homogeneous porous media is a first step toward developing efficient and accurate numerical methods for modelling unsaturated flows though heterogenous soils.

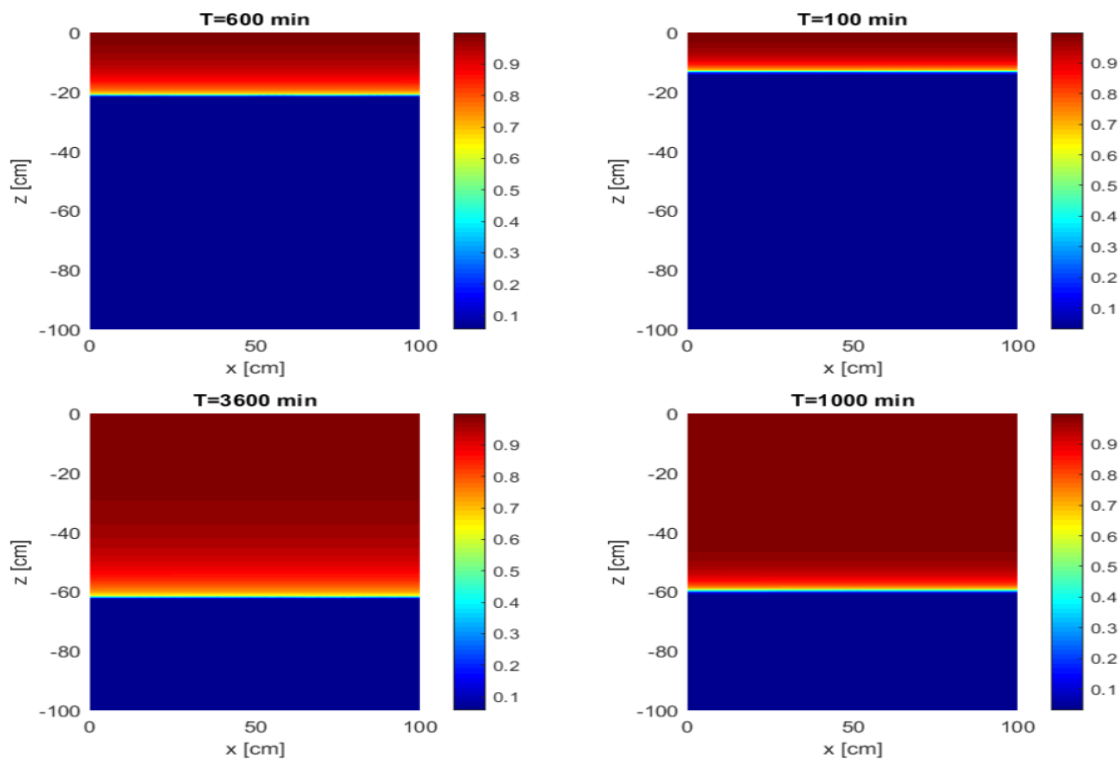

Figure 3: The time evolution of saturation for the sandy clay (left) and loam (right) soils.